\RequirePackage{fix-cm}
\documentclass[smallextended]{svjour3}       % onecolumn (second format)
\smartqed  % flush right qed marks, e.g. at end of proof
\usepackage{graphicx}
 \usepackage{mathptmx}      % use Times fonts if available on your TeX system
 \usepackage{algorithm, algorithmic}
\begin{document}

\title{Preconditioned Krylov subspace methods for sixth order compact approximations of the Helmholtz equation%\thanks{Grants or other notes
}
\author{Yury Gryazin
}
\institute{Y. Gryazin \at
              Mathematics Department, Idaho State University, 921 S 8th Ave, STOP 8085, Pocatello, Idaho 83209-8085, USA \\
              Tel.: 1(208)282-5293\\
              Fax: 1(208)282-2636\\
              \email{gryazin@isu.edu}           %  \\
}
\date{Received: date / Accepted: date}

\titlerunning{Preconditioned Krylov subspace methods for compact schemes.}
\maketitle

\begin{abstract}
In this paper, we consider an efficient iterative approach to the solution of the discrete Helmholtz equation with Dirichlet, Neumann and Sommerfeld-like boundary conditions based on a compact sixth order approximation scheme and preconditioned Krylov subspace methodology. A sixth order compact scheme for the 3D Helmholtz equation with different boundary conditions is developed to reduce approximation and pollution errors, thereby softening the point-per-wavelength constraint. The resulting systems of finite-difference equations are solved by different preconditioned Krylov subspace-based methods. In the majority of test problems, the preconditioned Generalized Minimal Residual (GMRES) method is the superior choice, but in the case of sufficiently fine grids a simple stationary two-level algorithm proposed in this paper in combination with a lower order approximation preconditioner presents an efficient alternative to the GMRES method. In the analysis of the lower order preconditioning developed here, we introduce the term  ``$k$-th order preconditioned matrix"  in addition to the commonly used ``an optimal preconditioner". The necessity of the new criterion is justified by the fact that the condition number of the preconditioned matrix $ AA^{-1}_p $ in some of our test problems improves with the decrease of the grid step size. In a simple 1D case, we are able to prove this analytically. This new parameter could serve as a guide in the construction of new preconditioners. The lower order direct preconditioner used in our algorithms is based on a combination of the separation of variables technique and Fast Fourier Transform (FFT) type methods. The resulting numerical methods allow efficient implementation on parallel computers. Numerical results confirm the high efficiency of the proposed iterative approach.

\keywords{Compact finite-difference schemes \and Krylov subspace methods \and GMRES \and FFT \and preconditioning}
\end{abstract}

\section{Introduction}
In recent years, the problem of increasing the resolution of existing numerical solvers has become an urgent task in many areas of science and engineering. Most of the existing efficient solvers for structured matrices were developed for lower-order approximations of partial differential equations. The need for improved accuracy of underlying algorithms leads to modified discretized systems and as a result to the modification of the numerical solvers (see e.g.\cite{zs}). 
	
The use of a lower order preconditioner for efficient implementation of high-order finite-difference and finite-element schemes has been under consideration for a long time (see e.g. \cite{hmmo}, \cite{orz}). In this paper, a compact sixth order approximation finite-difference scheme is developed, and a lower order approximation direct solver as a preconditioner for an efficient implementation of this compact scheme for the Helmholtz equation in the Krylov subspace method framework is considered. This approach allows us to utilize the existing lower order approximation solvers which significantly simplifies the implementation process of the higher resolution numerical methods. 

The model problem considered in the paper is the numerical solution of
the Helmholtz equation 
\begin{equation}
\nabla ^{2}u+k^{2}u=f , \ \  \textnormal{in}  \ \Omega, \label{problem}
\end{equation}
with the Dirichlet, Neumann, and/or Sommerfeld-like (radiation) boundary conditions 
\begin{eqnarray}
u=0,  \textnormal {on}  \ \partial \Omega_1, \nonumber  \label{bc} \\ 
u_{n}=0,  \textnormal {on} \  \partial \Omega_2, \\  
u_{n}-iku=0,  \textnormal {on} \ \partial \Omega_3 , \nonumber 
\end{eqnarray}
where $\Omega =\left\{ 0 \leq x,y,z\leq a\right\}$, $k$ 
is a complex valued constant, and $\partial \Omega_1,\partial \Omega_2$ and $\partial \Omega_3$ are different boundary sides of the rectangular computational domain $\Omega$.

It is known that for a given error level in the numerical approximation to the solution of the   Helmholtz equation, the quantity $(Re(k))^{p+1}h^{p}$ needs to be constant, where $  p $ is the order of finite-difference scheme and $h$ is the grid size. This phenomenon is known as ``pollution" \cite{bs,bgt2}. One way of reducing the pollution error is to increase the order of accuracy of the scheme. In this paper a sixth order compact finite-difference scheme is considered to address this problem. In the cases of the 3D Helmholtz and Dirichlet boundary conditions, this scheme was proposed by Sutmann in \cite{sut}. The 2D version of this scheme with Dirichlet and Neumann boundary conditions was developed in  \cite{nsd,st}. In this paper, the method is extended to the explicit compact sixth order approximation of Neumann and Sommerfeld-like boundary conditions in the 3D case. The extension of the known approach for the Neumann boundary conditions is straightforward. But in the case of the Sommerfeld-like boundary conditions the method is nontrivial and requires the introduction of a new auxiliary function.  In the case of the variable coefficient Helmholtz equation and Sommerfeld-type boundary conditions, the fourth order compact approximation scheme was considered in \cite{btt}. 

The resulting discretization leads to a system of linear
equations with block 27-diagonal structure. In general, the matrix of this system is neither positive definite nor Hermitian. Hence, most iterative methods either fail to
converge or converge too slowly, which is impractical. For the solution of this system we propose using a combination of Krylov subspace-based methods and the FFT preconditioner. Concerning other
approaches for the solution of the problem (\ref{problem})-(\ref{bc}), we refer to
Bayliss, Goldstein and Turkel  \cite{bgt1} for a preconditioned conjugate-gradient
algorithm, Kim \cite{k} for a domain decomposition method,
and Douglas, Hensley and Roberts \cite{dhr} for an ADI algorithm.
A very efficient multigrid method based on a shifted-Laplace preconditioner for the Helmholtz equation is presented in  \cite{umo}. 
The analysis of the multigrid algorithms in the case of nonsymmetric and indefinite problems can be found in \cite{bpx}.

On the other hand, the solution of this problem by a direct method
based on Gaussian elimination requires a prohibitive amount of additional
storage and computer time and thus has limited use. The most promising
results in the solution of a similar problem have been obtained by
preconditioned Krylov subspace methods in \cite{eg}. In this paper we
generalize some approaches developed for the second-order central difference discretization of the Helmholtz equation by Elman and O'Leary \cite{eo,eo1} and the author and others \cite{gkl} to the case of compact sixth order approximation scheme. 

The key to the fast convergence of the suggested iterative method is the choice of the preconditioning matrix. In our algorithm we use a preconditioner based on the second order central difference approximation of the Helmholtz equation and corresponding boundary conditions. The inversion of the preconditioning matrix at each step of the Krylov subspace method is done by a direct solver based on the FFT technique which requires $O(N^3\log{N})$ operations, where $N$ is the number of grid points in each direction. 

Numerical experiments with test problems demonstrate the high resolution of the sixth order compact finite-difference scheme, as well as the computational efficiency of our preconditioned Krylov subspace type numerical method. In most situations, the GMRES method demonstrates the best convergence properties. However, in the case of sufficiently fine grids, we propose using a simple stationary two-level method (see e.g.\cite{sn}).  This method was naturally constructed in the analysis of the convergence of the GMRES algorithm. So, the choice of the parameters in this method is not based on the minimization of the spectral radius of the iteration matrix on each step but rather on the construction of a particular linear combination from a Krylov subspace. This is why, for the purpose of this paper, we took the liberty of calling this approach the Simplified Krylov Subspace (SKS) method. On sufficiently fine grids this method requires much less processor time than the GMRES method, though the number of iterations until convergence is still larger than in the case of the GMRES method as should be expected. Also, since the implementation of SKS algorithm does not require the calculation of scalar products, it has greater potential than the GMRES method for implementation on parallel computers. A distinguishing feature of these approaches is that the number of iterations required for the convergence of Krylov subspace iterations decreases as the size of the discretized system increases. To explain this fact we introduce ``the order of the preconditioned matrix" as a parameter to quantify the rate of convergence of the condition number of the preconditioned matrix $ AA^{-1}_p $ to 1 on a sequence of grids. In some simple situation, we were able to find this parameter analytically. We believe that this parameter is more informative then the commonly used "an optimal preconditioner" (see e.g. \cite{l}, p.196) in the case of a lower order preconditioners and may be used as guide in the further development of preconditioners of similar type. But we must notice that even in this paper this parameter has limited application. 

This conclusion is based on the theoretical analysis of some simple situations and confirmed in our numerical experiments. Some preliminary numerical results on the application of this approach were presented at the 10th International Conference on Mathematical and Numerical Aspects of Waves, Vancouver, Canada, 24-29 July, 2011 \cite{g}.

The rest of the paper is organized as follows. In Section 2, the main idea of the proposed sixth order approximation compact method is presented in the case of the 1D problem. To analyze the convergence of the algorithms developed here, variations of known convergence estimates for the Krylov subspace methods are considered.  In this section, simplified approaches based on the Krylov subspace iterations are also presented.  Section 3 focuses on the development of  the sixth order compact approximation scheme in the case of Neumann and Sommerfeld-like boundary conditions. In Section 4, the effectiveness of the proposed algorithms is demonstrated on a series of test problems.

\section{A one-dimensional model problem} 
 
Let $A$ and $A_{p}$ be matrices derived from six and second order approximations to (\ref{problem}) and the Dirichlet boundary condition (\ref{bc}) using a mesh $x_{i}=ih, i=0,1,...,N+1, h=\frac{1}{N+1}$. The standard notation for the first and second order central differences at $ith$ grid point is given by
\begin{equation}
\delta_{x}u_i=\frac{u_{i+1}-u_{i-1}}{2h}, \ 
\delta_x^2u_i=\frac{u_{i+1}-2u_i+u_{i-1}}{h^2} \label{dif},
\end{equation}
where $u_i=u(x_i)$.
The difference operators $\delta_{y}$, $\delta_{z}$, $\delta_{y}^2$ and $\delta_{z}^2$ used in the following sections are defined similarly. 
The sixth order approximation to the second derivative at the $ith$ grid point can be written as 
\begin{equation} 
u''_i=\delta_x^2u_i - \frac{h^2}{12}u_i^{(4)} - \frac{h^4}{360}u_i^{(6)} + O(h^6). \label{T1D}
\end{equation} 
As usual, in the case of compact schemes, we need to use the original equation to find appropriate relations to eliminate the fourth and sixth derivatives in (\ref{T1D}). By using the second order central difference of the fourth derivative of $u$ at $x=x_i$ 
\begin{equation} 
\delta_x^2u_i^{(4)} = u^{(6)}_i  + O(h^6) 
\end{equation} 
and the expression for the fourth derivative of $u(x)$ from (\ref{problem})
\begin{equation} 
u^{(4)}_i=-k^2u''_i +f''_i, 
\end{equation} 
the relation (\ref{T1D}) can be expressed in the form
\begin{equation} 
\left (1- \frac{k^2h^2}{12}(1+ \frac{h^2}{30}\delta_x^2) \right )u''_i=\delta_x^2u_i -\frac{h^2}{12}\left(1+ \frac{h^2}{30}\delta_x^2\right)f''_i + O(h^6),  \label{T2D}
\end{equation} 
where $f_i=f(x_i)$.
After using (\ref{T2D}), the discretized system corresponding to the compact sixth order approximation scheme for (\ref{problem})-(\ref{bc}) can be presented as
\begin{eqnarray*}
&d_1U_{i-1}+d_2U_i+d_3U_{i+1}=F_{i},\  \  i=1,...,N, \\ 
&U_0 = 0, \ U_{N+1}=0,\\
&F_{i}=h^2 \left (1- \frac{7k^4h^4}{90} \right )f_i -\frac{k^4h^4}{360}(f_{i-1}+f_{i+1})+ \frac{7h^4}{90}f''_i+\frac{h^4}{360}(f''_{i-1}+f''_{i+1}), \\
&d_1=d_3=1-\frac{k^4h^4}{360}, d_2 = -2 + k^2h^2 -\frac{7k^4h^4}{90}, 
\end{eqnarray*}
where $U_i$ is the sixth order discrete approximation to $u(x_i)$. This system can also be rewritten in the form $AU=F$, where
 \begin{equation} 
A = d_1\Lambda+d_2I  \label{1DA} \\
\end{equation}  
 and 
\begin{equation} 
\Lambda =\left[\begin{array}{ccccc} 0&1&&\cdots&0\\
1&0&1&&\vdots\\ &\ddots&\ddots&\ddots&\label{LMB}\\ 
\vdots&&1&0&1\\
0&\cdots&&1&0\end{array}\right].\nonumber
\end{equation} 
In a similar way the preconditioning matrix based on the second order central difference approximation can be presented as
\begin{equation} 
A_p = \Lambda-(2 - k^2h^2 )I . \label{1DAP} \\
\end{equation} 
Finally the right preconditioned system can be written
\begin{eqnarray}
&AA_p^{-1}Y=F,  \label{AAP} \\
&A_pU=Y.  \nonumber
\end{eqnarray}

The main goal of this paper is to demonstrate the computational efficiency of a preconditioning technique based on the use of a lower order approximation discrete system in the numerical implementation of higher order compact finite-difference schemes. To consider this construction in more general form, we introduce the definition of a $kth$ order preconditioned system.

\begin{definition}
Given two nonsingular $N\times N$ matrices $A$ and $A_p$, a system in the form of (\ref{AAP})
 is said to be a $kth$ order preconditioned system if the $N\times N$ matrix $AA^{-1}_p$ is diagonalizable and can be expressed as $V^{-1}(I+h^kD)V$, where $h<1$; $D=$diag$(d_{11},...,d_{nn})$ with 
$ \displaystyle{ \max_i} |d_{ii}| < M$, where $M$ does not depend on $h$;
 and $V$ is an  $N\times N$ matrix of eigenvectors of $AA^{-1}_p$. 
\end{definition}

Using this definition we present the convergence analysis of the GMRES method applied to the system (\ref{AAP}) in the following theorem.
\begin{theorem} \label{Tm1}
Let a system in the form of (\ref{AAP}) be a $kth$ order preconditioned system. Then the $nth$ iteration $U^{(n)}$ of the GMRES method applied to this system satisfies the convergence estimate 
\begin{equation} 
\| r_n \|_2 \leq \kappa_2(V) (Mh^k)^n\| r_0\|_2, \label{est}
\end{equation} 
where $\kappa_2(V)=\| V^{-1}\|_2 \| V\|_2 $ and $r_n=F-AU^{(n)}$.
 \end{theorem} 

{\em Proof}. Since the matrix $AA^{-1}_p$ is diagonalizable, it is well known (see e.g.\cite{GV},\cite{s}) that the residual of $nth$ iteration of GMRES algorithm satisfies 
\begin{equation} 
\| r_n \|_2 \leq \kappa_2(V) \min_{p\in P_n, p(0)=1}\ \  \max_{i=1,...,N} |p(1+h^kd_{ii})|\cdot \| r_0\|_2. 
\end{equation} 
Consider the polynomial of $nth$ degree $\widehat{p}_n(x)=1-xy_{n-1}(x)$, where $y_{n-1}(x)=\displaystyle{ \sum_{k=0}^{n-1}}\alpha^n_{k+1}x^j$ and the coefficients $\alpha^n_{k}$ satisfy the equation  $\widehat{p}_n(1+x)=-\alpha^n_{n}x^n, n \ge 1$. It is easy to see that the coefficients $\alpha^n_{k}, \ k=1,...,n$ satisfy 
\begin{eqnarray}
&\displaystyle{ \sum_{k=1}^{n}\alpha^n_{k}}=1, \label{sys} \\
&\displaystyle{ \sum_{k=l-1}^{n}{k \choose {l-1} }\alpha^n_{k}}=0, \ l=2,...,n ,\nonumber
\end{eqnarray}
and the unique solution of this system is given by $ \alpha^n_k=(-1)^{k-1}{n \choose k }, k=1,...,n$. Then the convergence estimation for the GMRES method becomes
\begin{eqnarray} 
\| r_n \|_2 &\leq& \kappa_2(V) \min_{p\in P_n, p(0)=1} \max_{i=1,...,N} |p(1+h^kd_{ii})|\cdot \| r_0\|_2 \\ 
&\leq&  \kappa_2(V)  \max_{i=1,...,N} |\widehat{p}_n(1+h^kd_{ii})|\cdot \| r_0\|_2 \\
&\leq& \kappa_2(V) (Mh^k)^n\| r_0\|_2. \ \ \ \nonumber  
\end{eqnarray}
\endproof

The proof is somehow trivial and could be significantly simplified but we present it in such a form to use later in the construction of the simplified iterative method based on this proof.

In the same way we can derive another useful estimate expressed in the following corollary.
\begin{corollary}  \label{CTm1}
Let the matrix $AA^{-1}_p$ in the system (\ref{AAP}) be expressible as $V^{-1}(I+D)V$, where  $D=diag(d_{11},...,d_{nn})$, 
$M=\displaystyle{ \max_i} |d_{ii}| < 1 $ and $V$ is the  $N\times N$ matrix of eigenvectors of $AA^{-1}_p$.  Then the $nth$ iteration $U_n$ of the preconditioned GMRES method applied to this system satisfies the convergence estimate 
\begin{equation} 
\| r_n \|_2 \leq \kappa_2(V) M^n\| r_0\|_2. 
\end{equation} 
 \end{corollary} 
We can derive another well known estimate based on a standard application of Chebyshev polynomials (see e.g.\cite{s}). We consider only the case of real valued matrices. The following theorem provides the details of the proposed technique.
\begin{theorem} \label{Tm2}
Let the matrix $AA^{-1}_p$ in the system (\ref{AAP}) be  diagonalizable and expressible in the form $V^{-1}(I+D)V$, where  $D=diag(d_{11},...,d_{nn})$, $ -1< \widehat{m}=\displaystyle{ \min_i} \ d_{ii} \le \widehat{M}=\displaystyle{ \max_i} \ d_{ii}$ and $V$ is the  $N\times N$ matrix of eigenvectors of $AA^{-1}_p$.   Then the $nth$ iteration $U_n$ of the preconditioned GMRES method applied to this system satisfies the convergence estimate 
\begin{equation} 
\| r_n \|_2 \leq 2\kappa_2(V) \left( \frac{\widehat{M}-\widehat{m}}{4(1+\widehat{m})}\right ) ^n\| r_0\|_2. \label{CA}
\end{equation} 
 \end{theorem} 
 
{\em Proof}. This estimation immediately follows from the general complex valued matrix result (see e.g.\cite{s}, p.206). \endproof 

Next we will applied the developed convergence estimations to the analysis of the solution of the system (\ref{AAP}). 
\begin{theorem}\label{Tm3} 
Let $N\times N$ matrices $A$ and $A_p$ be defined by (\ref{1DA}) and (\ref{1DAP}) with $h<\frac{2\pi}{10k}$ and $\displaystyle{\min_{j=1,...,N} \left | \frac{ 4\sin^2(\frac{j\pi h}{2}) } {h^2}-k^2 \right |} \ge \delta_0 > 0$, where $ \delta_0$ is constant. Then $AA^{-1}_p$ is a second order preconditioned system and we have the convergence estimate.
\begin{equation} 
\| r_n \|_2 \leq \left ( \frac {k^4}{12\delta_0}h^2 \right ) ^n\| r_0 \|_2.  \label{1Dest} \\
\end{equation}
\end {theorem} 
{\em Proof}.
The matrices $A$ and $ A_p$ have the same set $\{V_j, j=1,...,N\}$ of eigenvectors, as the matrix $\Lambda$ from (\ref{LMB}). The eigenvectors are given by (see e.g.\cite{eo1})
\begin{equation} 
V_{j}^l=\sqrt{2h}\sin{(jl\pi h)}, j,l=1,...,N.
\end{equation} 
The matrix $V$ is unitary, so $\| V^{-1}\|_2= \| V\|_2 =1.$
The eigenvalues of the matrices $A_p$ and $A$ are given by
\begin{eqnarray} 
\lambda_{j}^p&=&-4\sin^2{\left (\frac{j\pi h}{2}\right)}+h^2k^2, \ \textnormal {for} \ j=1,...,N \ \ \textnormal {and} \\
\lambda_{j}&=&-4\sin^2{\left(\frac{j\pi h}{2}\right)}+h^2k^2-h^4k^4\frac{(14+\cos{(j\pi h)})}{180}, \ \textnormal {for} \  j=1,...,N \nonumber
\end{eqnarray} 
respectively. We can write the eigenvalues of the preconditioned matrix $AA_p^{-1}$ in the form 
\begin{eqnarray} 
	\frac {\lambda_{j}} { \lambda_{j}^p}=1+\frac { h^4k^4\frac{(14+\cos{(j\pi h)})}{180} }{4\sin^2{(\frac{j\pi h}{2})}-h^2k^2}, \ \textnormal {for} \ j=1,...,N. \label{1DEGV}
\end{eqnarray} 
The second term in the right hand side can be estimated using $ m \le | \frac {\lambda_{j}} { \lambda_{j}^p}-1| \le M $, where $m=13h^4k^4/720$ and $M=h^2k^4/(12\delta_0)$. So $AA^{-1}_p$ is a second order preconditioned matrix and using (\ref{est}) we obtain the estimate stated in the theorem.
\endproof 

The first condition of this theorem $ hk<2\pi/10 $ is a common natural restriction that comes to play even when one wants to visualize a numerical solution for qualitative analysis. It is just not an accurate representation of the solution if there are fewer then 5 points per half wave length. Moreover, to avoid the so-called the ``pollution" phenomenon Bayliss et. al. \cite{bgt2} introduced the restriction that $k^{p+1}h^p$ should be constant, i.e. to maintain a fixed accuracy of a scheme, the number of grid points must grow as $k^{1/p}$ where $p$ is the accuracy order of the scheme. 
So, if one satisfies the restriction that avoids ``pollution" phenomenon, the first condition of the theorem is satisfied almost automatically and it is not much of a restriction at all. 

The second condition of the theorem is a requirement that avoids the discrete spectrum of the operator. Since the resolvent set of a bounded linear operator is open, we can always do this by introducing a small change to $k^2$ if necessary. 

Next, we present how we can choose $\delta_0$ in some important cases. 

In the simplest case of $k^2 < 9$, we can take $\delta_0=1/3$. 

Consider how to get an estimate for $\delta_0$ in the case of a given $k > \pi$ and a sequence of grids with the grid sizes $ h_1>h_2>h_3\dots $ . The first restriction of the theorem  gives us the condition $kh <\frac{2\pi}{10}$. We can see that if the $j$ in the argument of $\displaystyle{\min_{j} \left | \frac{ 4\sin^2(\frac{j\pi h}{2}) } {h^2}-k^2 \right |}$ were to take on values of all real numbers from $1$ to $N$, then the minimum of the expression under consideration would be zero. 

Now let's denote $\alpha(j) = \frac{j\pi h}{2} $ and consider for which $\alpha$ the expression is zero if $kh =\frac{2\pi}{10}$. We consider the worst-case scenario to justify this approach for all $kh < \frac{2\pi}{10}$. It is trivial to find that $\alpha_0 = \sin^{-1}(\frac {kh}{2})=\sin^{-1}(\frac {\pi}{10})\approx .32 $. So, it is clear that the minimum of the function occurs either at $j_0=\lfloor  \frac{2\alpha_0}{\pi h} \rfloor $ or at $ j_1=j_0+1$ since $\sin(\alpha(j))$ is increasing on the interval $ 1<j<N$. There exists a neighborhood of $\alpha_0$ which includes both $\alpha(j_0)$ and $\alpha(j_1)$ such that $\alpha < .64$ which allows us to use the Taylor series to represent $\sin^2(\alpha)$ in this neighborhood. We use only first two terms in the series, i.e.
$ \sin^2(\alpha) = \alpha^2 - \frac{1}{3}\cos(2\xi) \alpha^4, 0< \xi < 0.64$. Now we can substitute this expression in the equation 
\begin{eqnarray}
\displaystyle{\min_{j=1,...,N} \left | \frac{ 4\sin^2(\frac{j\pi h}{2}) } {h^2}-k^2 \right |}= \min_{j=j_0,j_1} \left | j^2\pi^2 -
\frac{1}{12} \cos(2\xi)(j\pi)^4 h^2-k^2 \right | \nonumber \\
 \geq \min\left [ k^2 - j_0^2\pi^2, j_1^2\pi^2-k^2 - \frac{1}{12}(j\pi)^4 (h_1)^2 \right]  \geq \delta_0. \nonumber
\end{eqnarray}
So, we just need to avoid the values for $k$ and $h_1$ for which one of the terms in the brackets is nonpositive. In some situations we would need to choose a sufficiently small grid size step $h_1$ for the coarsest grid in a sequence. For example, in the case $k = 20$ and $h_1=1/32$,  $\min\left (400 - 6^2\pi^2, 7^2\pi^2-400 -\frac{1}{12} \frac{(7\pi)^4}{32^2} \right) > \min( 44,63)  = \delta_0$. This $\delta_0$ is independent of $h$ and can be used in our convergence analysis on a sequence of grids.

It follows from (\ref{1Dest}) that the number of iterations for the Krylov subspace algorithm decreases as the grid size $h$ decreases. This result proves that the algorithm developed here yields a very effective iterative technique for the implementation of a higher order approximation scheme. Later we will consider the extension of this method to 3D problems. 

The proof of Theorem \ref{Tm1} suggests a simplified version of the Krylov subspace method for the solution of (\ref{AAP}). Indeed, using the expression for the exact solution of the system (\ref{sys}), we have $\alpha^n_n=-\alpha^{n+1}_{n+1}$, for $n=1,...$ and the following recurrence expression for $y_n$ 
\begin{eqnarray} 
	y_n(z) =1-zy_{n-1}(z)+y_{n-1}(z), \textnormal {where} \ n \ge 1 .\nonumber
\end{eqnarray} 	
	This allow us to derive the SKS algorithm described in the table labelled Algorithm 1. 

\begin{algorithm}
\caption{Simplified Krylov Subspace (SKS) algorithm.  }
\begin{algorithmic} [1] \label{SKS}                  % enter the algorithmic environment
\STATE Let $Y^{(0)}=A_p^{-1}U^{(0)}=0$ be the initial approximation, let the tolerance be $tol$, and let the maximum number of iterations be $M$
\STATE $  r_0=F $, $Y =r_0$ and  $\delta_0= \| r_0 \|_2$
\STATE {$j=1$}
\WHILE{$j < M$ and $ \epsilon < tol$}
\STATE $U^{(j)}=A_p^{-1}Y$
\STATE $w=AU^{(j)}$
\STATE $r=r_0-w$
\STATE $\delta=\| r \|_2$
\STATE $\epsilon=\frac{\delta}{\delta_0}$
\STATE $Y=Y+r$
\STATE $j=j+1$
\ENDWHILE
\STATE  $U^{(j)}$ is the iterative solution of $AU=F$.
\end{algorithmic}
\end{algorithm} 
	
	This is a simple stationary two-level method (see e.g.\cite{zs}).  However, the choice of the parameters in this algorithm is not based on the minimization of the spectral radius of the iteration matrix on each step but rather on the construction of a particular linear combination of vectors from a Krylov subspace. This is why, for the purpose of this paper, we call this approach the Simplified Krylov Subspace (SKS) method. The proposed method does not require an orthogonalization procedure, i.e. requires no inner products, which is attractive for an implementation in a parallel computing environment. In addition, this method does not use estimates of maximum and minimum eigenvalues of the preconditioned matrix which is the drawback of the Chebyshev acceleration algorithm (see e.g. \cite{s}). The Chebyshev acceleration method is naturally constructed in the proof of the Theorem  \ref{Tm2} and we will compare the numerical effectiveness of the methods in Section 4. We should notice that in most numerical experiments, the GMRES method exhibits the best convergence properties, but in some situations the SKS algorithm proves to be more efficient. Some such examples are considered in Section 4. 
	
	The SKS algorithm also gives a good criterion for evaluating the quality of a preconditioner in the form of the order of preconditioned matrix (\ref{AAP}). In this paper we will calculate this as follows. Consider two grids with the grid size $h$ and $\gamma h,$ where $ \gamma >1$.  Let the $l_2-$norm of the residuals of the SKS method on first two iterations be $\| r^h_1 \|_2 $ and $\| r^h_2 \|_2 $ on the first grid and $\| r^{\gamma h}_1 \|_2 $ and $\| r^{\gamma h}_2 \|_2 $ on the second grid. Let $\epsilon^h=\frac {\| r^h_2 \|_2}{\| r^h_1 \|_2} $ and $\epsilon^{\gamma h}=\frac {\| r^{\gamma h}_2 \|_2}{\| r^{\gamma h}_1 \|_2} $. Now we can approximately calculate the order of preconditioned matrix by 
\begin{eqnarray} 	
\psi = \frac{\ln {\epsilon^{\gamma h} }  - \ln {\epsilon^{ h} }}{\ln { \gamma } }.  \label{PSI} 
\end{eqnarray} 		
	We consider this parameter in the discussion of the results of the numerical experiments.

\section{Three dimensional problems} 
\subsection{A sixth order approximation compact scheme} 
In this section we present a three dimensional compact sixth order approximation finite-difference scheme that was first introduced in \cite{sut}, and we develop a sixth order compact explicit approximation of the Neumann and Sommerfeld-like boundary conditions (\ref{bc}). In this discussion we consider a uniform grid   $ \Omega_h= \{ (x_i, y_j, z_k)| x_i=ih, y_j=jh, z=kh; i,j,k=0,...,N-1,  h=a/(N-1) \} $.

First, we consider a finite-difference approximation of (\ref{problem}) in the form 
\begin{equation} 
\delta_x^2u_{i,j,k} +\delta_y^2u_{i,j,k}+\delta_z^2u_{i,j,k}+ {k^2}u_{i,j,k} + T_{i,j,k}=f_{i,j,k} , 
\end{equation} 
where $u_{i,j,k}=u(x_i,y_j,z_k)$ and  $f_{i,j,k}=f(x_i,y_j,z_k)$ and 
\begin{eqnarray} 
 T_{i,j,k}=-\left[ \frac{h^2}{12}\left ( \frac{\partial^4u}{\partial x^4}+\frac{\partial^4u}{\partial y^4}+\frac{\partial^4u}{\partial z^4}\right ) + \frac{h^4}{360}\left ( \frac{\partial^6u}{\partial x^6}+\frac{\partial^6u}{\partial y^6}+\frac{\partial^6u}{\partial z^6} \right )\right]_{i,j,k} + O(h^6). \nonumber 
\end{eqnarray} 

Using the appropriate derivatives of (\ref{problem}) we can write the sixth order compact approximation of the Helmholtz equation in the form
\begin{eqnarray}
&&( \delta_x^2 + \delta_y^2 + \delta_z^2)U_{i,j,k}+\frac{h^2}{6}\left (1+ \frac{k^2h^2}{30}\right )( \delta_x^2\delta_y^2+ \delta_x^2\delta_z^2+ \delta_y^2\delta_z^2)U_{i,j,k}+ \nonumber\\
&&\frac{h^4}{30}\delta_x^2\delta_y^2\delta_z^2U_{i,j,k}+k^2\left (1- \frac{k^2h^2}{12}+\frac{k^4h^4}{360}\right )U_{i,j,k}= \label{scheme} \\
&&\left (1- \frac{k^2h^2}{12}+\frac{k^4h^4}{360} \right )f_{i,j,k}+\frac{h^2}{12}\left (1- \frac{k^2h^2}{30}\right )\nabla^2 f_{i,j,k}+\frac{h^4}{360}\nabla^4 f_{i,j,k} + 
\nonumber\\
&&\frac{h^4}{90}\left ( \frac{\partial^4f}{\partial x^2\partial y^2}+\frac{\partial^4f}{\partial x^2\partial z^2}+\frac{\partial^4f}{\partial y^2\partial z^2}\right )_{i,j,k} .\nonumber
\end{eqnarray}
For the two dimensional case this scheme was proposed in \cite{nsd}. In the three dimensional case with Dirichlet boundary conditions, it was developed in \cite{sut}. We will consider the iterative implementation of this scheme based on preconditioned Krylov subspace type algorithms.
\subsection{Boundary conditions} 
The implementation of the Dirichlet boundary conditions (\ref{bc}) is straightforward but the explicit compact approximation of  Neumann and Sommerfield-like boundary conditions requires careful consideration. The sixth order compact approximation of the Neumann boundary condition in the two dimensional case was considered in \cite{nsd}. Here we extend this approach to the three dimensional problem.  
\subsubsection{Neumann boundary conditions} 
In this subsection the Neumann boundary conditions are considered in the form 

\begin{equation} 
\left. \frac{\partial u}{\partial z} \right |_{z=0} = \beta(x,y). \label{3Dbcn}
\end{equation} 

We restrict our consideration to only one side of the computational domain $\Omega$.
By using Taylor series, we can derive the sixth order approximation formula
\begin{equation}
\delta_z u_{i,j,0} =\left.  \frac{ \partial u}{\partial z} \right |_{i,j,0}  + \frac{h^2}{6} \left.  \frac{ \partial^3 u}{\partial z^3} \right |_{i,j,0}  +  \frac{h^4}{120} \left.  \frac{ \partial^5 u}{\partial z^5} \right |_{i,j,0} +O(h^6). \label{T3DN}
\end{equation}
To express the third and fifth derivatives in (\ref{T3DN}) we can differentiate the original Helmholtz equation (\ref{problem}), assuming sufficient smoothness of the solution and the right hand side. After substitution of these derivatives into (\ref{T3DN}), the resulting expression is  
\begin{eqnarray}
&&\delta_z u_{i,j,0} =\left [  \frac{ \partial u}{\partial z} \right ]_{i,j,0}  + \frac{h^2}{6} \left [\frac{ \partial f}{\partial z}-\frac{ \partial^3 u}{\partial x^2\partial z}-\frac{ \partial^3 u}{\partial y^2\partial z}- k^2\frac{ \partial u}{\partial z} \right ]_{i,j,0}   \nonumber\\
&& + \frac{h^4}{120} \left [\frac{ \partial^3 f}{\partial z^3} -  
\frac{ \partial^5 u}{\partial z^3\partial x^2} - \frac{ \partial^5 u}{\partial z^3\partial y^2} - 
k^2\left ( \frac{ \partial f}{\partial z}-\frac{ \partial^3 u}{\partial x^2\partial z}-\frac{ \partial^3 u}{\partial y^2\partial z}- k^2\frac{ \partial u}{\partial z} \right ) \right  ]_{i,j,0} \label{DU}\\
 && +\ O(h^6).\nonumber
\end{eqnarray}
Next, we approximate the third order mixed derivatives in the second bracket of (\ref{DU}) by a fourth order approximation formula and in the second line using a second order approximation scheme, which yields
\begin{eqnarray}
\delta_z u_{i,j,0}& =&\left (1- \frac{ h^2k^2}{6}+ \frac{ h^4k^4}{120}  \right )\left [  \frac{ \partial u}{\partial z} \right ]_{i,j,0}  + 
\frac{ h^2}{6}\left (1-  \frac{ h^2k^2}{20}  \right )\left [  \frac{ \partial f}{\partial z} \right ]_{i,j,0} + \frac{ h^4}{120}\left [  \frac{ \partial^3 f}{\partial z^3} \right ]_{i,j,0}\nonumber\\
& - &\frac{h^2}{6} \left [\delta_z \delta_x^2 u +\delta_z \delta_y^2 u - \frac{h^2}{6}  \left (\frac{ \partial^5 u}{\partial x^2\partial z^3} + \frac{ \partial^5 u}{\partial y^2\partial z^3} \right )-\frac{h^2}{12}  \left ( \frac{ \partial^5 u}{\partial x^4\partial z}+\frac{ \partial^5 u}{\partial y^4\partial z} \right )\right ]_{i,j,0}   \nonumber\\
& + &\frac{h^4k^2}{120} \left [\delta_z \delta_x^2 u_{i,j,0} +\delta_z \delta_y^2 u_{i,j,0}\right ] - \frac{h^4}{120}\left [ \frac{ \partial^5 u}{\partial x^2\partial z^3} + \frac{ \partial^5 u}{\partial y^2\partial z^3}  \right ]_{i,j,0} +\ O(h^6).\nonumber
\end{eqnarray}
Now we can simplify this and write 
\begin{eqnarray}
\delta_z u_{i,j,0}&+&\frac{ h^2}{6}\left (1-  \frac{ h^2k^2}{20}  \right ) \left [  \delta_z \delta_x^2 u_{i,j,0} +\delta_z \delta_y^2 u_{i,j,0}  \right]  \nonumber \\
&=& \left (1- \frac{ h^2k^2}{6} + \frac{ h^4k^4}{120}  \right ) \left [  \frac{ \partial u}{\partial z} \right ]_{i,j,0} 
 + \frac{ h^2}{6}\left (1-  \frac{ h^2k^2}{20}  \right )\left [  \frac{ \partial f}{\partial z} \right ]_{i,j,0} + \frac{ h^4}{120}\left [  \frac{ \partial^3 f}{\partial z^3} \right ]_{i,j,0} \nonumber\\
& +& \frac{ 7h^4}{360}\left [ \frac{ \partial^3 f}{\partial z \partial^2 x} + \frac{ \partial^3 f}{\partial z \partial^2 y}\right ]_{i,j,0} 
- \frac{ 7h^4}{180} \left [  \frac{ \partial^5 u}{\partial x^2\partial y^2\partial z}  \right ]_{i,j,0}  \nonumber\\
&-& \frac{ 7h^4k^2}{360}\left [\frac{ \partial^3 u}{\partial x^2\partial z} + \frac{ \partial^3 u}{\partial y^2\partial z}  \right ]_{i,j,0} 
- \frac{h^4}{180} \left [\frac{ \partial^5 u}{\partial x^4\partial z} + \frac{ \partial^5 u}{\partial y^4\partial z}  \right ]_{i,j,0} +\ O(h^6).\nonumber
\end{eqnarray}
Finally, to complete the explicit scheme for the boundary conditions, we present the previous equation in the form
\begin{eqnarray}
\delta_z u_{i,j,0}&+&\frac{ h^2}{6}\left (1+  \frac{ h^2k^2}{30}  \right ) \left [  \delta_z \delta_x^2 u_{i,j,0} +\delta_z \delta_y^2 u_{i,j,0}  \right]  + \frac{ h^4}{30} \left [  \delta_z \delta_x^2 \delta_y^2 u_{i,j,0}  \right] = \nonumber\\ 
&=&\left (1- \frac{ h^2k^2}{6} + \frac{ h^4k^4}{120}  \right ) \left [  \frac{ \partial u}{\partial z} \right ]_{i,j,0} 
+ \frac{ h^2}{6}\left (1-  \frac{ h^2k^2}{20}  \right )\left [  \frac{ \partial f}{\partial z} \right ]_{i,j,0} + \frac{ h^4}{120}\left [  \frac{ \partial^3 f}{\partial z^3} \right ]_{i,j,0}  \nonumber\\
&+& \frac{ 7h^4}{360}\left [ \frac{ \partial^3 f}{\partial z \partial^2 x} + \frac{ \partial^3 f}{\partial z \partial^2 y}\right ]_{i,j,0} - \frac{ h^4}{180} \left [  \frac{ \partial^5 u}{\partial x^2\partial y^2\partial z}  \right ]_{i,j,0}  \nonumber\\
&-& \frac{ h^4k^2}{180}\left [\frac{ \partial^3 u}{\partial x^2\partial z} + \frac{ \partial^3 u}{\partial y^2\partial z}  \right ]_{i,j,0} - \frac{h^4}{180} \left [\frac{ \partial^5 u}{\partial x^4\partial z} + \frac{ \partial^5 u}{\partial y^2\partial z}  \right ]_{i,j,0} +\ O(h^6).\nonumber
\end{eqnarray}
Now by using (\ref{dif}) and replacing $u_{i,j,k}$ with $U_{i,j,k}$, we can write
\begin{eqnarray}
&& U_{i,j,-1}+\frac{ h^2}{6}\left (1+  \frac{ h^2k^2}{30}  \right ) \left [ \delta_x^2 U_{i,j,-1} + \delta_y^2 U_{i,j,-1}  \right]  + \frac{ h^4}{30} \delta_x^2 \delta_y^2 U_{i,j,-1}    \label{3DN} \\  
&=&U_{i,j,1}+\frac{ h^2}{6}\left (1+  \frac{ h^2k^2}{30}  \right ) \left [ \delta_x^2 U_{i,j,1} + \delta_y^2 U_{i,j,1}  \right]  + \frac{ h^4}{30} \delta_x^2 \delta_y^2 U_{i,j,1}   \nonumber\\
&-&2h \left (1- \frac{ h^2k^2}{6} + \frac{ h^4k^4}{120}  \right ) \beta_{i,j,0}- 
\frac{ h^3}{3}\left (1-  \frac{ h^2k^2}{20}  \right )\left [  \frac{ \partial f}{\partial z} \right ]_{i,j,0} - \frac{ h^5}{60}\left [  \frac{ \partial^3 f}{\partial z^3} \right ]_{i,j,0}  \nonumber\\
&+&\frac{ 7h^5}{180}\left [ \frac{ \partial^3 f}{\partial z \partial^2 x} + \frac{ \partial^3 f}{\partial z \partial^2 y}\right ]_{i,j,0} +  
\frac{ h^5}{90} \left [  \frac{ \partial^4 \beta}{\partial x^2\partial y^2}    + k^2\left (\frac{ \partial^2 \beta}{\partial x^2} + \frac{ \partial^2 \beta}{\partial y^2}  \right )  
+ \frac{ \partial^4 \beta}{\partial x^4} + \frac{ \partial^4 \beta}{\partial y^4}  \right ]_{i,j,0}.  \nonumber
\end{eqnarray}
This formula allows us to eliminate the term $U_{i,j,-1}$ in (\ref{scheme}). The explicit implementation of the Neumann boundary conditions on the other boundary pieces can be conducted in a similar way. 
\subsubsection{Sommerfield-like boundary conditions}
In this section we consider the implementation of the Sommerfeld-like boundary conditions  (\ref{bc}) which are often used in scattering problems 
\begin{equation} 
\left. \frac{\partial u}{\partial z} +iku\right |_{z=0} = 0.\label{3DS} \\
\end{equation} 
The difficulty in using similar methods to those described in the previous section can be seen in Equation (\ref{3DN}). The last term in this equation requires calculation of the fourth derivative of the right hand side of (\ref{3Dbcn})   and if it depends on an unknown variable $u$, the use of a compact sixth order approximation scheme becomes problematic. To avoid such a difficulty,
we introduce a new variable $v=e^{ikz}u$. Then the Helmholtz equation becomes
\begin{equation}
\frac{\partial^2 v}{\partial x^2}+\frac{\partial^2 v}{\partial y^2}+\frac{\partial^2 v}{\partial z^2}-2ik\frac{\partial v}{\partial z}=\bar{f} , \label{3DMH} \\
\end{equation}
where $ \bar{f} =e^{ikz}f$. Now (\ref{3DS}) can be rewritten in the form
\begin{equation} 
\left. \frac{\partial v}{\partial z} \right |_{z=0} = 0.\label{3DSM} \\
\end{equation} 
The sixth order approximation of (\ref{3DSM}) becomes
\begin{equation}
\delta_z v_{i,j,0} =\left.  \frac{ \partial v}{\partial z} \right |_{i,j,0}  + \frac{h^2}{6} \left.  \frac{ \partial^3 v}{\partial z^3} \right |_{i,j,0}  +  \frac{h^4}{120} \left.  \frac{ \partial^5 v}{\partial z^5} \right |_{i,j,0} +O(h^6). \label{3DT} \\
\end{equation}
Assuming sufficient smoothness of $v$ and $\bar{f}$, and using (\ref{3DMH}) and (\ref{3DSM}), we can express the third and fifth derivatives of $v$ in the form 

\begin{eqnarray}
\left. \frac{\partial^3 v}{\partial z^3} \right |_{z=0} &=& \left[ \frac{\partial \bar{f} }{\partial z} +2ik \frac{\partial^2 v}{\partial z^2 } \right ]_{z=0} \label{3D3DR} \\
\left. \frac{\partial^5 v}{\partial z^5} \right |_{z=0} &=& -4ik\left[ 2k^2\frac{\partial^2 v }{\partial z^2} +\frac{\partial^4 v}{\partial z^2 \partial x^2} +\frac{\partial^4 v}{\partial z^2 \partial y^2} \right ]_{z=0}   \label{3D5DR} \\
&-& \left[ 4k^2\frac{\partial \bar{f} }{\partial z} -2ik\frac{\partial^2 \bar{f}}{\partial z^2 } +\frac{\partial^3 \bar{f}}{\partial z \partial x^2}+\frac{\partial^3 \bar{f}}{\partial z \partial y^2}-\frac{\partial^3 \bar{f}}{ \partial z^3} \right ]_{z=0} .  \nonumber 
\nonumber \end{eqnarray}
To preserve the sixth order compact approximation in (\ref{3DT}) we need to use the fourth order approximation for  $ \frac{\partial^2 v}{\partial z^2 } $ in (\ref{3D3DR}). The second order approximation is sufficient for the derivatives of $v$ in (\ref{3D5DR}). First, let's consider the fourth finite-difference compact approximation for  $ \frac{\partial^2 v}{\partial z^2 } $. It can be presented in the form
\begin{eqnarray}
&& \left. \frac{\partial^2 v}{\partial z^2} \right |_{i,j,0} = \delta_z^2 v_{i,j,0} -\frac{ h^2}{12}\left.  \frac{ \partial^4 v}{\partial z^4} \right |_{i,j,0} + O(h^4) \label{3DT1} \\
&=&\delta_z^2 v_{i,j,0} +\frac{ h^2}{12}\left[ 4k^2 \frac{ \partial^2 v}{\partial z^2}+  \frac{\partial^4 v}{\partial z^2 \partial x^2} +\frac{\partial^4 v}{\partial z^2 \partial y^2} -2ik\frac{\partial \bar{f} }{\partial z} -\frac{\partial^2 \bar{f}}{\partial z^2 }   \right ]_{i,j,0} + O(h^4)
 \nonumber \\
 &=&\delta_z^2 v_{i,j,0} +\frac{ h^2}{12}\left[ 4k^2 \delta_z^2 v_{i,j,0}+ \delta_z^2  \delta_x^2 v_{i,j,0} + \delta_z^2  \delta_y^2 v_{i,j,0}  -\left (2ik\frac{\partial \bar{f} }{\partial z} +\frac{\partial^2 \bar{f}}{\partial z^2 }   \right )_{i,j,0}\right ] + O(h^4).  
\nonumber \end{eqnarray}
Now we can use (\ref{3DSM}), (\ref{3DT1}) and the second order approximation of (\ref{3D5DR}) to express (\ref{3DT}) in the form 
\begin{eqnarray}
&& \delta_z v_{i,j,0} -\frac{ ikh^2}{3}\left( 1+\frac{ 2k^2h^2}{15} \right ) \delta_z^2 v_{i,j,0}  + \frac{ ikh^4}{180}\left(  \delta_z^2  \delta_x^2 v_{i,j,0} + \delta_z^2  \delta_y^2 v_{i,j,0} \right )  =\bar{F}_{i,j,0}, \label{3DIMP}\\
&&\bar{F}_{i,j,0} = \left( \frac{ h^2}{6}+\frac{ k^2h^4}{45} \right )  \left. \frac{\partial \bar{f} }{\partial z} \right |_{i,j,0}- \frac{ ikh^4}{90} \left.  \frac{\partial^2 \bar{f}}{\partial z^2 } \right |_{i,j,0} + \frac{ h^4}{120}\left[\frac{\partial^3 \bar{f}}{ \partial z^3}- \frac{\partial^3 \bar{f}}{\partial z \partial x^2}-\frac{\partial^3 \bar{f}}{\partial z \partial y^2} \right ]_{i,j,0} . 
\nonumber \end{eqnarray}
This equation gives an implicit compact sixth order approximation for (\ref{3DSM}). Next, the explicit implementation of the boundary conditions (\ref{3DSM}) similar to (\ref{3DN}) is developed by using the equation 
\begin{eqnarray}
\ \ \ \ \ \ \ 0&=&h^2 \mu_1\left [ \frac{\partial^3 v}{\partial z \partial x^2}+\frac{\partial^3 v}{\partial z \partial y^2} \right ]_{i,j,0}   \label{3DT2} \\
&=&  \mu_1h^2\left[\delta_z \delta_x^2 v_{i,j,0} +\delta_z \delta_y^2 v_{i,j,0}- \frac{ h^2}{6}\left ( \frac{\partial^5 v}{\partial z^3 \partial x^2} +\frac{\partial^5 v}{\partial^3 z \partial y^2}\right )_{i,j,0} \right ]+ O(h^6)\nonumber \\
&=& \mu_1h^2 \left[\delta_z \delta_x^2 v_{i,j,0} +\delta_z \delta_y^2 v_{i,j,0}\right ] -\frac{ik\mu_1 h^4}{3}\left [ \frac{\partial^4 v}{\partial z^2 \partial x^2} +\frac{\partial^4 v}{\partial^2 z \partial y^2}\right ]_{i,j,0}   \nonumber \\
&+&\frac{\mu_1 h^4}{6}\left [2\mu_2\frac{\partial^5 v}{ \partial x^2\partial y^2\partial z} -\frac{\partial^3\bar{f} }{\partial z \partial x^2}-\frac{\partial^3\bar{f} }{\partial z \partial y^2}\right ]_{i,j,0} + O(h^6) . 
\nonumber \end{eqnarray}
Here, $\mu_1$ and $\mu_2$ are parameters. We also use  that $\left. \frac{\partial^5 v}{\partial z \partial x^2 \partial y^2}\right |_{z=0}=\left. \frac{\partial^5 v}{\partial z \partial x^4 }\right |_{z=0}=0$ and $\left. \frac{\partial^5 v}{\partial z \partial y^4 } \right |_{z=0}=0$. We leave the first of these derivatives in the right hand side but drop the latter two derivatives. Now, by using the second order finite-difference approximation in (\ref{3DT2}) and adding it to the (\ref{3DIMP}), we derive the formula
\begin{eqnarray}
&& \mu_3 \left [ \delta_z v_{i,j,0} -\frac{ ikh^2}{3}\left( 1+\frac{ 2k^2h^2}{15} \right ) \delta_z^2 v_{i,j,0} \right ]  + \mu_3 \frac{ ikh^4}{180}\left[  \delta_z^2  \delta_x^2 v_{i,j,0} + \delta_z^2  \delta_x^2 v_{i,j,0} \right ] \\ 
&&+\mu_1 \mu_3 \left[  \delta_z  \delta_x^2 v_{i,j,0} + \delta_z  \delta_x^2 v_{i,j,0}  -  \frac{ ikh^4}{3}\left(  \delta_z^2  \delta_x^2 v_{i,j,0} + \delta_z^2  \delta_x^2 v_{i,j,0} \right )  \right ] \nonumber  \\
&&+\mu_1 \mu_2 \mu_3 \delta_z  \delta_x^2 \delta_z v_{i,j,0} = \mu_3\bar{F} + \mu_1 \mu_3 \frac{ h^4}{6}\left[ \frac{\partial^3 \bar{f}}{\partial z \partial x^2}+\frac{\partial^3 \bar{f}}{\partial z \partial y^2} \right ]_{i,j,0} . 
\nonumber \end{eqnarray}
Next we choose parameters $ \mu_1, \mu_2$ and $\mu_3 $ to match the coefficients in (\ref{scheme}). This choice is determined by the conditions 
\begin{eqnarray}
&& \mu_3=\frac {1}{ 1 + \frac{ 2ikh}{3}\left( 1+\frac{ 2k^2h^2}{15} \right )} , \ \mu_1=\frac {\frac{ ikh}{90} + \frac{ 1}{6\mu_3}\left( 1+\frac{ k^2h^2}{30} \right )}{ 1 + \frac{ 2ikh}{3}} , \ \mu_2 = \frac{ 1}{10\mu_1\mu_3}.
 \end{eqnarray}
Using these parameters and replacing $v_{i,j,l}$ with $U_{i,j,l}e^{ikz_l}$, we write
\begin{eqnarray}
\ \ \ \ \ &&U_{i,j,-1}+\frac{ h^2}{6}\left (1+  \frac{ h^2k^2}{30}  \right ) \left ( \delta_x^2 + \delta_y^2  \right)U_{i,j,-1}   + \frac{ h^4}{30} \delta_x^2 \delta_y^2 U_{i,j,-1}  \label{3DDISCS}\\ 
&=&e^{2ikh} \mu_3 \left [ h^2\left (\mu_1+ \frac{ 2ikh}{3}\left ( \frac{ 1}{60}-\mu_1 \right ) \right) \left ( \delta_x^2 + \delta_y^2  \right)U_{i,j,1}  \right. \nonumber \\
&&\ \ \ \ \ \ \ \ \ \ \ \ + \frac{ h^4}{30\mu_3} \left.  \delta_x^2 \delta_y^2 U_{i,j,1}  +\left (1 - \frac{ 2ikh}{3} -  \frac{ 4ih^3k^3}{45}  \right )U_{i,j,1}\right]   \nonumber \\
&+&e^{ikh} \mu_3 \left [\frac{ 4ikh}{3} \left (1 +  \frac{ 2h^2k^2}{15}  \right )U_{i,j,0}- \frac{ 4ikh^3}{3}\left ( \frac{ 1}{60}-\mu_1 \right ) \left ( \delta_x^2 + \delta_y^2  \right)U_{i,j,0}  \right]  \nonumber \\
&-&e^{ikh} \mu_3 \left [2h\bar{F}_{i,j,0} + \mu_1 \frac{ h^5}{3}\left( \frac{\partial^3 \bar{f}}{\partial z \partial x^2}+\frac{\partial^3 \bar{f}}{\partial z \partial y^2} \right )_{i,j,0} \right ]. \nonumber 
\end{eqnarray}
This equation provides an explicit compact sixth-order approximation for the boundary condition (\ref{3DS}). At the upper boundary $ z=a $, the Sommerfeld-like boundary condition can be written $ \left. \frac{\partial u}{\partial z} - iku\right |_{z=a} = 0 $ and the auxiliary substitution becomes $v=e^{-ikz}u$. The rest of the derivation is very similar to the derivation for the case of the lower boundary. Implementations of the sixth order compact approximation of the boundary conditions in the other directions are also similar to the calculations already presented.  

\subsection{Fast Fourier Preconditioner}

The main idea of this paper is to utilize the existing lower-order approximation direct solvers as preconditioners in the implementation of a higher resolution scheme. One of the most popular methods for the approximate solution of the Helmholtz equation is the second-order central difference scheme. We consider the cases with Dirichlet or Neumann boundary conditions on the sides of the computational domain $ \Omega $, i.e. $ u(0,y,z)=u(a,y,z)=$ $u(x,0,z)=u(x,a,z)=0 $ or $  \left. \frac {\partial u}{\partial x}\right |_{x=0}=\left. \frac {\partial u}{\partial x}\right |_{x=a}=\left. \frac {\partial u}{\partial y}\right |_{y=0}=\left. \frac {\partial u}{\partial y}\right |_{y=a} =0$. At the top $z=a$ and the bottom $z=0$ of the computational domain, we assign one of the three boundary conditions under consideration: Dirichlet, Neumann or Sommerfeld-like boundary conditions (\ref{bc}). In the case of Neumann and Sommerfeld-like boundary conditions, we use the second order central finite-difference approximation of the first derivative on all boundaries. To express the preconditioner in general form we introduce the $N\times N$ matrices
\begin{equation} 
\Lambda_{N}^\beta =\left[\begin{array}{ccccc} 0&1+\beta&&\cdots&0\\
1&0&1&&\vdots\\ &\ddots&\ddots&\ddots&\\ 
\vdots&&1&0&1\\
0&\cdots&&1+\beta&0\end{array}\right]\nonumber
\end{equation} 
and $ D_N^\alpha=diag(\alpha, 0, ... , 0, \alpha) $, where $\alpha, \beta$ and $N$ are parameters depending on boundary conditions. Next let $ A_{\beta,N}^0 =  \Lambda_N^\beta - 2I$ and $A_{\alpha,\beta,N} =  \Lambda_N^\beta - (2-k^2h^2)I+ D_N^\alpha$. The preconditioning matrix can now be written in the form 
\begin{equation}
 A_p=A_{\beta,N}^0\otimes I \otimes I +I \otimes A_{\beta,N}^0 \otimes I + I \otimes I \otimes A_{\alpha_z,\beta_z,N_z}, \label{3DAP} \\
\end{equation} 
where $\beta=0$ in the case of Dirichlet boundary conditions at the sides of the computational domain and $ \beta = 1$ in the case of Neumann boundary conditions on the same boundaries. The number of grid points in the $z$-direction $N_z$ is chosen in such a way that the grid step size is the same in all three directions. 
This condition is imposed so that the lower order preconditioner and the sixth-order approximation compact scheme developed in the previous section use the same grid points. Depending on the boundary conditions at the top and the bottom of the computational domain, the parameters $\alpha_z$ and $\beta_z$ are defined in the following way: in the case of Dirichlet boundary conditions $\alpha_z=0$ and $\beta_z=0$; in the case of Neumann boundary condition $\alpha_z=0$ and $\beta_z=1$; and in the case of the Sommerfeld-like boundary conditions $\alpha_z=2ikh$ and $\beta_z=1$. At each step of the iterative process, the solution to the second equation in (\ref{AAP}) with the matrix (\ref{3DAP}) can be obtained by using FFT-type algorithms in $O(N^2N_z\log_2{N})$ operations.
\subsection{Convergence analysis of 3D algorithm} 

In this section, the convergence of the proposed algorithms is considered in the case of the 3D Helmholtz equation with the Dirichlet boundary conditions (\ref{problem}), (\ref{bc}) on the rectangular computational domain $ \Omega =\left\{ 0\leq x,y,z\leq 1\right\}$ with uniform grid size $h=\frac {1}{N+1}$. To insure the uniqueness of the solution to the original boundary value problem, we assume that the coefficient $k^2$ is in the resolvent set of the Laplace operator defined on an appropriate function space. We also impose a common restriction on the number of grid points per wave length, i.e. $ kh < \frac {2\pi}{10} $. The convergence analysis of the proposed algorithms follows the same lines as in Theorem $\ref{Tm3}$. The preconditioning matrix is given by $(\ref{3DAP})$ with $ \beta =\alpha_z = \beta_z=0$ and $N_z=N$. To express the right hand side of the sixth order approximation scheme given in (\ref{scheme}), we will use the notation introduced in (\ref{3DAP}). In addition, we define the $N\times N$ matrix $ A_{\beta,N}^1 =  \Lambda_N^\beta - 4I$. Then the resulting $N^3\times N^3$ matrix $A$ can be written as 
\begin{eqnarray}
&&A=A_{0,N}^0\otimes I \otimes I +I \otimes A_{0,N}^0 \otimes I + I \otimes I \otimes A_{0,N}^0+ \\ 
&& \frac{1}{6}\left [1+  \frac{ h^2k^2}{30}  \right ](\Lambda_N^0 \otimes A_{0,N}^1 \otimes I +I \otimes A_{0,N}^1 \otimes A_{0,N}^1 -4 I \otimes I \otimes I)+\nonumber \\
&& \frac{1}{30}A_{0,N}^1 \otimes A_{0,N}^1 \otimes A_{0,N}^1+ \left [h^2k^2-\frac{ h^4k^4}{12}+\frac{ h^6k^6}{360}   \right ]  I \otimes I \otimes I.
\nonumber 
\end{eqnarray}
Both matrices $A$ and $A_p$ have the same set of orthonormal eigenvectors  $V_{i,j,l}$ (see e.g.\cite{eo1}):
\begin{eqnarray} 
V_{i,j,l}^{m,n,s}&=& \left( 2h \right)^{3/2}\sin{(im\pi h)}\sin{(jn\pi h)}\sin{(ls\pi h)}, \\
&& 1 \leq i,j,l,m,n,s \leq N. \nonumber
\end{eqnarray} 
The matrix of eigenvectors is unitary so the $l_2$-norm of $V$ and $V^{-1}$ is one. To present the eigenvalues of the matrices, we introduce the notation $s_r=\sin^2{(\frac{r\pi h}{2})}, r=i,j,l$. Then the eigenvalues of the matrices $A_p$ and $A$ in the 3D case are given by
\begin{eqnarray} 
\lambda_{i,j,l}^p&=&-4\left [ s_i^2 +s_j^2+s_l^2\right ]+h^2k^2, \textrm{and} \\
\lambda_{i,j,l}&=&-4\left [ s_i^2 +s_j^2+s_l^2 \right ]+h^2k^2
+ \frac{8}{3}\left [1+  \frac{ h^2k^2}{30}  \right ] \left [ s_i^2s_j^2 + s_i^2s_l^2+s_l^2 s_j^2\right ]\nonumber \\
&-& \frac{32}{15}s_i^2 s_j^2 s_l^2 - \frac{ h^4k^4}{12}+  \frac{ h^6k^6}{360}, i,j,l=1,...,N,\textrm{respectively}.\nonumber
\end{eqnarray} 
Unfortunately, the eigenvalues of these matrices do not satisfy the hypotheses of Theorem \ref{Tm1}, but we still can apply Corollary \ref{CTm1} to prove the convergence of the GMRES method and the SKS algorithm by showing that $ \lambda_{i,j,l}/\lambda_{i,j,l}^p=1+d_{i,j,l}$, where $|d_{i,j,l}|<1$, for all $i,j,l=1,...,N$ in some simple cases. 

Let's assume that $N\ge 1$, i.e. $h \le \frac{1}{2}$ and $k^2<25$ . Then $s_r^2=\sin^2(\frac{\pi rh}{2})\ge \sin^2(\frac{\pi h}{2}) \ge 2h^2$, for $r=i,j,l$. It is also easy to see  that $ \frac{8}{3} s_i^2s_j^2 \ge \frac{32}{15}s_i^2 s_j^2 s_l^2 $ and $ \frac{8}{3} (s_i^2s_l^2+s_j^2s_l^2) \ge \frac{64h^4}{3} \ge \frac{k^4h^4}{12}$ for the given values of $k$. 
Now the eigenvalues of the preconditioned matrix $AA_p^{-1}$ can be written in the form
\begin{eqnarray} 
\frac {\lambda_{i,j,l}}{\lambda_{i,j,l}^p}&=&1+ \left [ \frac{8}{3}\left (1+  \frac{ h^2k^2}{30}  \right ) \left [ s_i^2s_j^2 + s_i^2s_l^2+s_l^2 s_j^2\right ] -\frac{32}{15}s_i^2 s_j^2 s_l^2 - \frac{ h^4k^4}{12}+ \frac{ h^6k^6}{360}\right ]\ \ \ \  \label{3DEGV}\\
&/&\lambda_{i,j,l}^p,\  i,j,l=1,...,N.\nonumber
\end{eqnarray} 
This gives us a presentation of the eigenvalues of the preconditioned matrix in the form $ \lambda_{i,j,l}/\lambda_{i,j,l}^p=1+d_{i,j,l}$, where $|d_{i,j,l}|<3/4$, for all $i,j,l=1,...,N$. Now we can apply Corollary \ref{CTm1} and obtain an estimate for the convergence rate of convergence of the GMRES and SKS algorithms : 
\begin{equation} 
\| r_n \|_2 \leq (3/4)^n\| r_0 \|_2 , \label{3Dest}
\end{equation} 
where  $r_n=F-AU^{(n)}$.

Similarly to the 1D case, it is possible to show that this estimate holds for $k > 6$ and a sequence of grids with the grid sizes $ h_1>h_2>h_3\dots $ if 
\begin{eqnarray} 
\displaystyle{\min_{i,j,l=1,...,N} \left | \frac{ 4\left ( s_i^2 +s_j^2+s_l^2\right ) } {h^2}-k^2 \right |} \ge \delta_0 > 0
\end{eqnarray} 
for some $ \delta_0=const $ and sufficiently small $h_1$. As in the case of the 1D problem if follows from the fact that the resolvent set of a bounded linear operator is open. The proof is somewhat tedious but elementary, so we skip it. In the next table we present lower and upper bounds $m$ and $M$ for $ \lambda_{i,j,l}/\lambda_{i,j,l}^p-1$ for different values of $k$ and for different grid step sizes. Also, the last two rows of the table  present  possible choices for $h_1$ and $\delta_0$ for different $k$.

\begin{center}
\textbf{Table 1. 3-D Dirichlet problem eigenvalue bounds(different $k^2$)}
\end{center}

\begin{center}
\begin{tabular}{|r|r|r|r|r|r|}
\hline
       & $k=10$ & $k=20$ & $ k=30$ & $k=40$ & $k=50$\\  \hline
$h=1/64$ & -.015 \ \ $\vert$ 0.49 &-1.22 $\vert$ .49     &-1.02 $\vert$ .88&-34.7 $\vert$   \ 180  & -10.2 $\vert$   55.0   \\ \hline
$h=1/128$ & -.0038\  $\vert$ 0.49 & -.158 $\vert$ .49 &-.140 $\vert$ .64& -2.81 $\vert$ 0.90 & - 2.21 $\vert$  1.74  \\ \hline
$h=1/256$ &-.00094 $\vert$ 0.49& -.035 $\vert$ .49   &-.031 $\vert$ .49& -0.22 $\vert$ 0.69 & -11.7    $\vert$\ \ \  .56 \\ \hline
$h=1/512$ &  -.00023 $\vert$ 0.49& -.009 $\vert$ .49 &-.008 $\vert$ .49& -.048 $\vert$ 0.49 & -.298    $\vert$ \ \  .49 \\ \hline
$\delta_0 \ \ \ \ $ & 8.00 & 4.00  & 1.80 & 1.12 & 2.9 \\ \hline
$h_1 \ \ \ \ $ & 1/64 & 1/128  & 1/256 & 1/256 &  1/512 \\ \hline
\end{tabular}
\end{center}

This result is weaker then in 1D case, i.e. it says that convergence is independent of the grid step, but it does not improve with the decrease of $h$. In this case, $A_p$ is said to be ``an optimal preconditioner" (see e.g. \cite{l}, p.196). This situation is typical for multigrid-type algorithms but it is not what we would expect based on $1D$ example. In the following section we present the results of numerical experiments in the case of several test problems which confirm that actual convergence of the presented methods accelerate with the decrease of the grid size, i.e. the convergence in the numerical experiments significantly exceeds the estimate (\ref{3Dest}). 

\section{Numerical Results}

In this section we present the results of numerical experiments which demonstrate the quality of the proposed numerical methods. These algorithms were implemented in Matlab 7.11.0 on an iMac with an Intel Core i7, 2.93 GHz  processor. We also use the standard programing implementation of the GMRES method in Matlab (gmres function). 
\subsection{1D test problem}
In the first  series of numerical experiments, we consider the convergence of the compact sixth order approximation algorithm on a sequence of grids. 
In these tests we focus on the numerical solution of the 1D Dirichlet problem for the Helmholtz equation with zero boundary conditions on the interval $ 0\le x \le 1$ and $k=20$.
The computational grid in our experiments is defined by $x_i=ih, i=0,...,N+1$ where $h=\frac{1}{N+1}$ and $N$ is the number of grid points. We consider the function $u(x_i)= x_i(1-x_i)\cos(k\pi x_i), i=1,...,N$ as the exact solution of the original boundary value problem with the right hand side $f(x_i)=-(2+k^2(\pi^2-x_i(1-x_i))\cos(k\pi x_i)+ 2k\pi(2x_i-1)\sin(k\pi x_i), i=1,...,N$. 
The first column of Table 2 displays the step size $h$ of the grids in our experiments. Columns 2-4 of the table present the number of iterations until convergence required by the preconditioned GMRES method, the SKS method (Algorithm \ref{SKS}) and the Chebyshev acceleration (CA) algorithm \cite{s}. We use the stoppage criterion  $ ||F-AU^{(n)}||_2/||F||_2\le 10^{-10}$, where $U^{(n)}$ is the iterative solution on $nth$ iteration. Column 5 reports the approximation error of the numerical solution $Err_6=\max_{1 \le j \le N}|U^{(n)}_j-u(x_j)|$. In Column 6, the second order approximation error $Err_2$ of the numerical solution obtained by using the preconditioner as a solver is presented for comparison. The number of grid points is chosen  to satisfy the requirement $h<2\pi/10k$. The values of the minimal and maximum eigenvalues of the preconditioned matrix used by the Chebyshev acceleration algorithm are calculated exactly by (\ref{1DEGV}). 

\begin{center}
\textbf{Table 2. 1-D Dirichlet problem}
\end{center}

\begin{center}
\begin{tabular}{|r|r|r|r|r|r|}
\hline
 Grid Size ({\em h})  & GMRES & SKS & CA & $Err_6$ & $Err_2$ \\  \hline\hline
$1/32$ & 6&$8$ & $6$ &$3.88*10^{-03}$  &$1.16*10^{-01}$ \\ \hline
$1/64$ & 4& $4$ &$5$ & $5.08*10^{-05}$  &$2.4*10^{-02}$ \\ \hline
$1/128$ & 3&$3$ & $4$ &$7.61*10^{-7}$  &$5.7*10^{-03}$ \\ \hline
$1/258$ & 3&$3$ & $3$ &$1.14*10^{-8}$  &$1.4*10^{-03}$ \\ \hline
\end{tabular}
\end{center}

These experiments demonstrate the sixth order convergence of the numerical solution to the exact solution of the original boundary value problem. Theorem \ref{Tm3} suggests that the residual on $nth$ iteration in the presented iterative algorithms decreases with the increase in the number of grid points. The results of the numerical experiments confirm this conclusion. We have already shown that the order of the preconditioned matrix in this case is two. The numerical value of the constant $\psi$ calculated by using (\ref{PSI}) is $1.99$. The 1D numerical experiments confirm the effectiveness of the iterative approaches under consideration and allows us to expect that the same properties hold for 3D problems. 

\subsection{3D test problems}
In the next series of experiments we consider different boundary value problems for  the 3D Helmholtz equation (\ref{problem}),(\ref{bc}). The computational domain is given by $\Omega =\left\{ (x,y,z)| 0 \leq x,y,z\leq 1\right\}$ with a uniform grid in all three directions $x_i=ih, i=1,...,N_x,$ $ y_j=jh, j=1,...,N_y$ and $z_k=kh, k=1,...,N_z$. The numbers of grid points $N_x, N_y$ and $N_z$ in $x,y$ and $z$-directions are chosen to provide uniform grid size. 
\subsubsection{Dirichlet problem}
First, we consider the Helmholtz equation with boundary conditions $u=0$ and $k=20$. This test is similar to the 1D case with $N_x=N_y=N_z=N$ and $h=1/(N+1)$. As the solution of the boundary value problem we choose the function  $ u(x,y,z)= \phi_1(x)\phi_2( y)\phi_3(z)$, where $\phi_1( x) =x^3(1-x)^3 , \phi_2( y) =y(1-y)\cos(k\pi y), \phi( z) = \sin(k\pi z),$ and $ 0 \leq x,y,z\leq 1$.

In the first series of tests we consider the convergence of the compact high order approximation finite-difference scheme (\ref{scheme}) on a sequence of grids as well as convergence properties of the GMRES, SKS and CA iterative algorithms.  The iterative process is stopped when the \textit{initial} residual is reduced by a factor of $10^{-10}$. The values of the minimal and maximumal eigenvalues of the preconditioned matrix used by the Chebyshev acceleration algorithm are calculated exactly using (\ref{3DEGV}). The description of the columns of Table 3 is similar to the description of Table 2.

\begin{center}
\textbf{Table 3. 3-D Dirichlet problem}
\end{center}

\begin{center}
\begin{tabular}{|r|r|r|r|r|r|}
\hline
 Grid Size ({\em h})  & GMRES & SKS & CA & $Err_6$ & $Err_2$ \\  \hline\hline
$1/64$ & 5&$10$ & $23$ &$3.03*10^{-04}$  &$4.47*10^{-06}$ \\ \hline
$1/128$ & 3& $6$ &$15$ & $8.13*10^{-05}$  &$6.35*10^{-08}$ \\ \hline
$1/256$ & 3& $ 5$ & $14$ &$2.06*10^{-05}$  &$9.68*10^{-10}$ \\ \hline
\end{tabular}
\end{center}
The test results confirm the sixth order approximation of the compact finite-difference scheme. Though, in our theoretical analysis (\ref{3DEGV}) we could not prove the desired $kth$ order preconditioning property in the case of the 3D Dirichlet problem, the SKS algorithm still exhibits acceleration of the convergence with the decrease of grid size. We also mention that the total computer time required for convergence for the SKS and GMRES methods is approximately the same and the Chebyshev acceleration algorithm requires about four times more computer time until convergence than the other two iterative algorithms do.  
The SKS algorithm has also greater potential for efficient implementation on parallel computers than the GMRES method does. 
 
One potential application of the proposed methods is to electromagnetic scattering problems. In these problems, the typical value of the coefficient $k$ is  in the range between $10$ and $50$. For example, the propagation of an electromagnetic wave with 1 $GHz$ frequency in a vacuum is described by the Helmholtz equation with $k^2 \approx 20.9 \ m^{-1} $. In the next table we present the convergence of the algorithms presented here on the sequence of grids with different values $k$. In every cell of the table we present the number of iteration for three methods under consideration : GMRES, SKS, and CA. The stoppage criterion is the same as in the previous series of experiments. If any of these three methods fails to converge in 100 iterations, we indicate this by using symbol ``$>100$" and in the case of divergence we use ``div" symbol. 

\begin{center}
\textbf{Table 4. 3-D Dirichlet problem (different $k^2$)}
\end{center}

\begin{center}
\begin{tabular}{|r|r|r|r|r|r|}
\hline
       & $k=10$ & $k=20$ & $ k=30$ & $k=40$ & $k=50$\\  \hline
$h=1/64$ & 4$\vert$ 7 $\vert$ 14 &5$\vert$ 10 $\vert$ 23 & 6$\vert$ 13 $\vert$ 25 &12$\vert$ div $\vert$ $>100$  &19$\vert$ div $\vert$ $>100$   \\ \hline
$h=1/128$ & 3$\vert$ 5 $\vert$ 14 &3$\vert$ \ \ 6 $\vert$ 15 & 4$\vert$ \ \ 8 $\vert$ 17 &4$\vert$ \ 10 $\vert$ \ \ \ \ \ \  69 &5$\vert$ \ 11 $\vert$  $>100$  \\ \hline
$h=1/256$ & 3$\vert$ 4 $\vert$ 12& 3$\vert$ \ \ 5 $\vert$ 14 & 3$\vert$ \ \ 6 $\vert$ 13 & 3$\vert$ \ \ \ 6 $\vert$\  \ \ \ \ \ \ 22 &4$\vert$ \ \ \ 7 $\vert$  $\ \ \ \ \ \ 63$ \\ \hline
$\psi \ \ \ \ $ & 2.02 & 2.00  & 1.98 & 2.00 & 2.02 \\ \hline
\end{tabular}
\end{center}
From these numerical experiments we can see that the preconditioned GMRES and SKS methods demonstrate excellent convergence properties. We also observe that the convergence of the proposed algorithms improves with the increase in the number of the grid points. In the analysis of the algorithm's convergence (\ref{3Dest}), it was shown that the number of iterations does not increase with the increase of the grid points. But the last row in the table indicates that the parameter $ \psi $ (\ref{PSI}) for all values of $k$ is approximately 2, which indicates that similarly to the 1D case the convergence of the 3D algorithm not only does not depend on the grid size but accelerates with the increase of the number of grid points.  This also indicates that even in the case when the SKS method does not converge, one can expect convergence on a finer grid. In all numerical tests, the parameter $\psi$ is calculated by using the first two iterations of the SKS method on the grids with the grid sizes $h=\frac{1}{128}$ and $h=\frac{1}{256}$.

From the results presented, it is clear that slow convergence and difficulties arising in the choice of the parameters $\widehat{m}$ and $\widehat{M}$ in the Chebyshev acceleration algorithm (\ref{CA}) make this method a poor choice for the implementation of the developed high-order approximation scheme. These numerical tests illustrate the fact that the Chebyshev polynomial is not optimal on the discrete spectrum and this has a dramatic effect on the iterative methods based on this polynomial. So, in the next numerical experiments we focus on the convergence properties of the GMRES and SKS methods in our numerical framework. 

\subsubsection{Dirichlet-Neumann boundary conditions}
In the next series of experiments we replace the Dirichlet boundary conditions $ u=0$ in (\ref{bc}) with $ u=u_b$ at the top of the computational domain ($z=1$) and with the Neumann boundary conditions $\frac{ \partial u}{\partial z}=0$ at $ z=0$. On the other boundaries of the rectangular computational domain $\Omega =\left\{ 0 \leq x,y,z\leq 1\right\}$ we still use the Dirichlet boundary conditions $u=0$. As a test function that satisfies the boundary conditions we consider  $ u(x,y,z)= \phi_1(x)\phi_2( y)\phi_3(z)$, where $\phi_1( x) =x^3(1-x)^3 , \phi_2( y) =y(1-y)\cos(k\pi y), \phi( z) = \cos(k\pi z),$ and $ 0 \leq x,y,z\leq 1$. Also, $ u_b=\phi_1(x)\phi_2( y)\phi_3(1), 0 \leq x,y\leq 1$. The number of grid points is chosen to provide the same grid step in all three directions, i. e. $N_x=N_y=N_z-1$. The stoppage  criterion is the same as in the previous numerical experiments. Table 5 provides the results of the numerical tests on a sequence of grids with the coefficient $k=20$. As already mentioned we focus only on two best algorithms: GMRES and SKS. The rest of Table 5 is similar to Table 3.

\begin{center}
\textbf{Table 5. 3-D Dirichlet-Neumann problem}
\end{center}

\begin{center}
\begin{tabular}{|r|r|r|r|r|r|}
\hline
Grid Size ({\em h})  & GMRES & SKS & $Err_6$ & $Err_2$ \\  \hline\hline
$1/64$ & 8&$11$ & $4.65*10^{-06}$  &$3.67*10^{-04}$ \\ \hline
$1/128$ & 5& $8$ & $6.61*10^{-08}$  &$8.74*10^{-05}$ \\ \hline
$1/256$ & 4&$7$ & $1.01*10^{-09}$  &$2.17*10^{-05}$ \\ \hline
\end{tabular}
\end{center}
The results of this series of numerical tests indicate the sixth order convergence of the resulting compact finite-difference scheme including explicit approximation of the Neumann boundary conditions.  As expected, the number of iterations decreases with the increase of the number of grid points in both the GMRES and SKS approaches. In the next table we present the number of iterations required for convergence for both algorithms on a sequence of grids and with different $k$ coefficients in the Helmholtz equation. The content of Table 6 is similar to the content of Table 4, except for the number of iterative approaches under consideration. We present only the number of iterations for the GMRES and SKS methods. 

\begin{center}
\textbf{Table 6. 3-D Dirichlet-Neumann problem (different $k^2$)}
\end{center}

\begin{center}
\begin{tabular}{|r|r|r|r|r|r|}
\hline
       & $k=10$ & $k=20$ & $ k=30$ & $k=40$ & $k=50$\\  \hline
$h=1/64$ & \ \ 6\ \ \ $\vert$ \ \ 10 &8\ \ \ $\vert$ \ \ 11 & 18\ \ \ $\vert$ div &36\ \ $\vert$ \ \ div &$> 600$\ \ $\vert$ \ \ div   \\ \hline
$h=1/128$ &4\ \ \ $\vert$ \ \ \ 9 &5\ \ \ $\vert$ \ \ \ 8 & 8\ \ \ $\vert$ \ \ \ $9$ &14\ \ $\vert$ \ \ div  &28\ \ \ $\vert$ \ \ div  \\ \hline
$h=1/256$ & 3\ \ \ $\vert$ \ \ \  8& 4\ \ \ $\vert$ \ \ \ 7 &5\ \ \ $\vert$  div  & 9\ \ $\vert$ \ \  div   &7\ \ \ $\vert$ \ \ div  \\ \hline
$h=1/512$ & 2\ \ \ $\vert$ \ \ \  6& 3\ \ \ $\vert$ \ \ \ 6 & 3\ \ \ $\vert$\ \ \ \ 6  & 5\ \ $\vert$ \ \ \ \ 6   \ &4\ \ \ $\vert$ \ \ \ \ 6  \\ \hline
$\psi \ \ \ \ $ & 2.57 & 2.18  & 2.39 & 1.65  & 1.94 \\ \hline
\end{tabular}
\end{center}
It follows from this table that the GMRES approach is much more efficient on coarse grids but on finer grids, the processor time for the SKS method is the same or even smaller than for the GMRES algorithm. For example, on the grid with $h=1/512$ and $k=50$ the corresponding processor times for the GMRES and SKS methods are  $1467 \ sec$ and $1052 \ sec$. Also, as we mentioned before, the SKS method has much greater potential than the GMRES algorithm for efficient implementation on parallel computers. 
We also expect that even if the SKS method does not converge on coarser grids, by reducing the grid size one can achieve convergence of this algorithm. We can see this in the example of the convergence of the methods under consideration in the case of $k \ge 30$. The values of the parameter $\psi$ vary significantly  for different values of $k$ but $\psi>0$ still indicates that the convergence of the iterative algorithms improves with the increase in the number of the grid points. 

\subsubsection{Dirichlet-Sommerfeld-like boundary conditions}

In the last series of numerical experiments, we consider the boundary value problem for the Helmholtz equation (\ref{problem}) with a combination of Sommerfeld-like(radiation) boundary conditions (\ref{bc}) at $z=0$ and $z=1$, and Dirichlet boundary conditions at all other boundaries of the rectangular computational domain $\Omega =\left\{ 0 \leq x,y,z\leq 1\right\} $.
The source function $f$ in (\ref{problem}) selected that the true
solution is $u( x,y,z) =\phi_1(x)\phi_2( y)\phi_3(z)$, where $\phi_1( x) =x^3(1-x)^3 , \phi_2( y) =y(1-y)\cos(k\pi y),$$ \phi_3( z) =e^{ik(z+1)}+e^{-ik(z-1)}-2$ and  $0 \le x,y,z \le 1.$
Note that the analytic solution satisfies the radiation boundary condition
(\ref{bc}). The numbers of grid points in $x,y, z-$directions are chosen such that $N_x=N_y=N_z-2$. First we consider the convergence of our sixth order approximation scheme on the sequence of grids with $k=20$. The main goal of this series of numerical experiments is to investigate the convergence properties of the new explicit compact sixth order approximation of the Sommerfeld-like boundary conditions  (\ref{3DDISCS}) proposed in this paper. As in previous test runs, the iterative process is stopped when the $l_2-$norm of the \textit{initial}
residual is reduced by a factor of $10^{-10}$. The results are presented in Table 7. The description of data presented in this table is the same as in the case of Table 5.

\begin{center}
\textbf{Table 7. 3-D Dirichlet-Sommerfeld problem}
\end{center}

\begin{center}
\begin{tabular}{|r|r|r|r|r|r|}
\hline
Grid Size ({\em h})  & GMRES & SKS & $Err_6$ & $Err_2$ \\  \hline\hline
$1/64$ & $6$&$9$ & $1.25*10^{-06}$  &$1.35*10^{-03}$ \\ \hline
$1/128$ & $4$& $7$ & $1.84*10^{-08}$  &$3.23*10^{-04}$ \\ \hline
$1/256$ & $3$ &$5$ & $2.82*10^{-10}$  &$7.99*10^{-05}$ \\ \hline
\end{tabular}
\end{center}
In this series of numerical experiments, both algorithms exhibit the same convergence properties as in the previous series of tests, i.e. the number of iterations decreases as the number of grid points increases. From the table, we can also see the sixth order convergence of the approximate solution to the exact solution of the boundary value problem on a sequence of grids. These results confirm the sixth order approximation of the compact explicit scheme proposed for the numerical implementation of the Sommerfeld-like boundary conditions.

In the last table, we present the number of iterations required for the convergence of the GMRES and SKS methods for different coefficients $k$ in the Helmholtz equation with the same boundary conditions used in the previous series of numerical tests. Data presented in the next table are similar to the data in Table 6. 

\begin{center}
\textbf{Table 8. 3-D Dirichlet-Sommerfeld problem (different $k^2$)}
\end{center}

\begin{center}
\begin{tabular}{|r|r|r|r|r|r|r|}
\hline
       & $k=10$ & $k=20$ & $ k=30$ & $k=40$ & $k=50$& $k=35.7+0.43i $ \\  \hline
$h=1/64 \ $& 5 $\vert$ 7& 6 $\vert$ 9& 11 $\vert$ 43& 16 $\vert$ div& 119 $\vert$ div& 14 $\vert$ \ \ \ \ \ \ 54\\ \hline
$h=1/128$ & 4 $\vert$ 5& 4 $\vert$ 7& 7 $\vert$ \  8& 7 $\vert$ div& 12 $\vert$ 11& 7 $\vert$ \ \ \ \ \ \ \ 9\\ \hline
$h=1/256$ & 3 $\vert$ 4& 3 $\vert$ 5& 5 $\vert$ \ 6& 5 $\vert$ \ \ 6& 6 $\vert$ 7& 5 $\vert$ \ \ \ \ \ \ \ 6\\ \hline
$h=1/512$ & 2 $\vert$ 4& 3 $\vert$ 5& 4 $\vert$ \ 5& 3 $\vert$ \ \ 5& 4 $\vert$ 6& 4 $\vert$ \ \ \ \ \ \ \ 5\\ \hline
$\psi \ \ \ \ $ & 15 & 13  & 13 & 15  & 11 & 13 \ \ \ \ \  \ \ \\ \hline

\end{tabular}
\end{center}
In the last column of Table 8, we present the results of calculations when $k$ is a complex valued coefficient. The value of this coefficient approximately corresponds to the case of propagation of electromagnetic waves with the $1\ GHz$ frequency in dry soil. These results confirm the effective implementation of the compact sixth order approximation scheme by using the proposed Krylov subspace-type algorithms in the framework with the FFT based low order preconditioners. The last row of the table suggests that the acceleration of the convergence with the increase of the number of grid point is much stronger than the  acceleration in the previous series of experiments. It seems strongly dependent on the boundary conditions used in the numerical tests. In all our experiments, there is a clear connection between the parameter $\psi$ and the improvement of convergence with the decrease of the step side. But the usefulness of this parameter in the analysis of the quality of a preconditioner requires further analysis since in majority of our numerical experiments the preconditioner does not satisfy the condition in the Definition 1.  

The series of test problems considered suggests that in the majority of situations the preconditioned GMRES method is the most efficient choice for an effective implementation of the compact sixth order approximation scheme on the coarse grids but in the case of finer grids the SKS method in combination with lower order approximation preconditioner presents an efficient alternative to the GMRES method. This alternative could become even more valuable when the GMRES method experiences stagnation.

\section{Conclusions}
New 3D compact sixth order explicit finite-difference schemes for the approximation of Neumann and Sommerfeld-like boundary conditions on rectangular computational domains with uniform grid size were developed and implemented. Together with the compact sixth order approximation scheme for the Helmholtz equation proposed in \cite{sut}, these algorithms represent highly accurate methods for the solution of boundary value problems for Helmholtz equations.

A new rapid iterative method based on preconditioned Krylov subspace methodology was developed for the implementation of the proposed compact finite-difference schemes. The strategy is based on a combination of higher order approximation schemes and a lower order approximation preconditioner. The analysis of some typical test problems reveals the attractive properties of the developed methods such as the decrease of the number of iteration until convergence with the increase of the number of the grid points or the size of the resulting matrix. This approach is especially attractive in situations in which the lower approximation solver already exists and the original boundary value problem calls for more accurate approximation.  

The typical time to produce the sixth order accuracy solution of the 3D Helmholtz equation with a combination of the Dirichlet and Sommerfeld boundary conditions
on a $512^3$ grid was just 30 minutes on a iMac using only one 2.93 GHz Intel Core i7 processor. The method was tested for realistic parameter ranges typical for electromagnetic scattering problems. We must notice that the Sommerfeld-like boundary is just a first order approximation of the Sommerfeld conditions on the unbounded domain. So, direct application of the higher order approximation for the Sommerfeld-like boundary condition to the solution of a scattering problem is not always justified. However, a straightforward extension of the approximation approach presented in this paper could be applied for approximation of the absorbing boundary conditions (see e.g.  \cite{umo}) or can be used in the implementation of the  perfectly matched layer (PML) boundary conditions (see e.g.  \cite{hrt}).

\section*{Acknowledgments}
The author thanks Professor T. Payne and Professor D. Stowe for valuable comments and discussions.

\bibliographystyle{spmpsci}      % mathematics and physical sciences


\end{document}